\documentclass[11pt,twoside]{amsart}
\usepackage[bookmarksnumbered, plainpages]{hyperref}
\usepackage{amsthm}
\usepackage{amsmath}
\usepackage{amsfonts}
\usepackage{amssymb}
\usepackage{float}
\usepackage[dvipsnames]{xcolor}

\newtheorem{thm}{Theorem}[section]
\newtheorem{cor}[thm]{Corollary}

\theoremstyle{definition}

\theoremstyle{remark}

\begin{document}

\title{Normal Subgyrogroups of Certain Gyrogroups}

\author{\bf S. Mahdavi, A. R. Ashrafi and M. A. Salahshour$^{\star}$}

\thanks{$^{\star}$ Corresponding author (Email: salahshour@iausk.ac.ir)}

\address{\textbf{Soheila Mahdavi and Ali Reza Ashrafi:} Department of Pure Mathematics, Faculty of Mathematical Sciences, University of Kashan, Kashan 87317$-$53153, I. R. Iran}

\address{\textbf{Mohammad Ali Salahshour:} Department of Mathematics, Savadkooh Branch, Islamic Azad University, Savadkooh, I. R. Iran}

\begin{abstract}
Suppose that $(T, \star)$ is a groupoid  with a left identity such that each element $a \in T$ has a left inverse. Then $T$ is called a \textit{gyrogroup} if and only if $(i)$ there exists a function $gyr: T \times T \longrightarrow Aut(T)$ such that for all $a, b, c \in T$, $a \star (b \star c) =  (a \star b) \star gyr[a,b]c$, where $gyr[a,b]c = gyr(a,b)(c)$; and $(ii)$ for all $a, b \in T$, $gyr[a,b]$ $=$ $gyr[a \star b, b]$. In this paper, the structure of normal subgyrogroups of certain gyrogroups are investigated.

\vskip 3mm


\vskip 3mm

\end{abstract}

\maketitle

\section{Introduction}
Gyrogroup theory started in 1988  by Ungar \cite{6} in which he proved that the set of all 3-dimensional relativistically admissible velocities possesses a group-like structure in which the group-like operation is given by the standard relativistic velocity composition law.  In another paper \cite{7}, he has shown that the Thomas rotation, in turn, gives rise to a non-associative group-like structure for the set of relativistically admissible velocities. Nowadays this non-associative group-like structure is known as a gyrogroup. In an algebraic language,  if  $(T, \star)$ is a groupoid  with a left identity such that each element $a \in T$ has a left inverse. Then $T$ is called a \textit{gyrogroup} if and only if $(i)$ there exists a function $gyr: T \times T \longrightarrow Aut(T)$ such that for all $a, b, c \in T$, $a \star (b \star c) =  (a \star b) \star gyr[a,b]c$, where $gyr[a,b]c = gyr(a,b)(c)$; and $(ii)$ for all $a, b \in T$, $gyr[a,b]$ $=$ $gyr[a \star b, b]$.  It is easy to see that it is a generalization of a group by defining the gyroautomorphisms to be the identity automorphism.

Let $T$ be a gyrogroup and let $H$ be a non-empty subset of $T$. If $H$ is a group under the induced operation of $T$, then $H$ is called a \textit{subgroup} of $T$, and if $H$ is a gyrogroup under the induced operation of $T$, then we use the term \textit{subgyrogroup} for $H$. If $H$ is the kernel of a homomorphism from $T$ to another gyrogroup, then $H$ is called a \textit{normal subgyrogroup} of $T$; see \cite{2} for more details.

Suppose that $(H^+,\oplus)$ is a gyrogroup, $H^-$ is a non-empty set disjoint from $H^+$ such that $|H^+| = |H^-|$ and $\varphi: H^+ \longrightarrow H^-$ is bijective. Set $G = H^+ \cup H^-$ and $a^- = \varphi(a^+)$. For arbitrary elements $a^\varepsilon, b^\delta \in G$, we define:

\[ a^{\epsilon}\otimes  b^{\delta}=
 \begin{cases}
  a^{+} \oplus  b^{+}    &  \epsilon = \delta =+ \ or  \ \epsilon = \delta = -  \\
(a^{+} \oplus b^{+})^{-}   &  \epsilon = + , \delta = - \ or \ \epsilon = - ,  \delta = +\\
\end{cases}.
\]

\noindent The function $gyr_G: G \times G \longrightarrow Aut(G)$ is given by

\[gyr_{G}[a^{\epsilon},b^{\delta} ](t^{\gamma} )=
\begin{cases}
gyr_{H^{+}}[a^{+},b^{+}](t^{+})  & \gamma =+ \\
(gyr_{H^{+}}[a^{+},b^{+}](t^{+})) ^{-} & \gamma =-
\end{cases},
\]
where $a^{\epsilon}, b^{\delta}$ and $t^{\gamma}$ are arbitrary elements of $G$.

The aim of this paper is to prove the following theorem:

\begin{thm}\label{t1.1}
A non-empty subset $M$ of $G$ is a normal subgyrogroup if and only if one of the following  conditions are satisfied:
\begin{enumerate}
\item $M \unlhd H^+$;

\item there exists $N^+ \unlhd H^+$ and $L^- \subseteq H^-$ such that $M = N^+ \cup L^-$ and for each $x, y \in L^-$, $x \otimes y \in N^+$. Also, $N^+ \cap L^+ = \emptyset$ and $N^+ \cup L^+ \unlhd H^+$, where $L^+ = \varphi^{-1}(L^-)$;

\item $M = N^+ \cup N^-$ such that $N^+ \unlhd H^+$ and $N^- = \varphi(N^+)$.
\end{enumerate}
\end{thm}

Throughout this paper, our notations are standard and can be taken mainly from the books \cite{8,9}. We refer the readers to consult the survey article \cite{5} for a complete history of gyrogroups.

\section{Preliminary Results}
The following result of Suksumran \cite[Theorem 32]{2} is crucial throughout this paper:

\begin{thm}\label{t2.1}
Let $(T,\star)$ be a gyrogroup containing a subgyrogroup $H$. Then $H$ is normal in $T$ if and only if for all $a, b \in T$, $a \star (H \star b)=(a \star b) \star H=(a \star H) \star b$.
\end{thm}

The present authors \cite{1} obtained the structure of the subgyrogroups of $G$ which is important in finding its  normal subgyrogroups.

\begin{thm}\label{t2.2}
With above notations, $(G,\otimes)$ is a gyrogroup. A non-empty subset $B$ of $G$ is a subgyrogroup of $(G,\otimes)$ if and only if one of the following three conditions hold:
\begin{enumerate}
\item[(a)] $B \leq H^+$;

\item[(b)] there exists $A^+ \leq H^+$ and $L^- \subseteq H^-$ such that $B = A^+ \cup L^-$ and for each $x, y \in L^-$, $x \otimes y \in A^+$. Moreover, $A^+ \cap L^+ = \emptyset$ and $A^+ \cup L^+ \leq H^+$, where $L^+ = \varphi^{-1}(L^-)$;

\item[(c)] $B = A^+ \cup A^-$ such that $A^+ \leq H^+$ and  $A^- = \varphi(A^+)$.
\end{enumerate}
\end{thm}

\begin{cor}\label{c2.3}
The gyrogroup $(G,\otimes)$ satisfies the following conditions:
\begin{enumerate}
  \item $H^{+} \unlhd  G$.
  \item If $N^+ \unlhd H^+$, then $N^+ \unlhd G$
\end{enumerate}
\end{cor}

\begin{proof} The case (1) is \cite[Theorem 2.10]{1}. To prove (2), we assume that $a^{\epsilon},b^{\delta}$ are arbitrary elements of $G$. We consider two cases as follows:
\begin{enumerate}
\item[(i)] \textit{$\epsilon=\delta=+$ \ or \ $\epsilon=\delta=-$}. By definition,
\begin{eqnarray*}
  (a^{+}\otimes b^{+})\otimes N^{+} &=& (a^{+} \oplus b^{+})\otimes N^{+}=(a^{+} \oplus b^{+})\oplus N^{+} \\
  (a^{+}\otimes N^{+})\otimes b^{+} &=& (a^{+}\oplus N^{+})\otimes b^{+}=(a^{+}\oplus N^{+})\oplus b^{+} \\
  a^{+}\otimes (N^{+}\otimes b^{+}) &=& a^{+}\otimes(N^{+}\oplus b^{+})=a^{+}\oplus(N^{+}\oplus b^{+}) \\
  (a^{-}\otimes b^{-})\otimes N^{+} &=& (a^{+} \oplus b^{+})\otimes N^{+}=(a^{+} \oplus b^{+})\oplus N^{+} \\
  (a^{-}\otimes N^{+})\otimes b^{-} &=& (a^{+}\oplus N^{+})^-\otimes b^{-}=(a^{+}\oplus N^{+})\oplus b^{+} \\
  a^{-}\otimes (N^{+}\otimes b^{-}) &=& a^{-}\otimes(N^{+}\oplus b^{+})^-=a^{+}\oplus(N^{+}\oplus b^{+})
\end{eqnarray*}
Note that by our assumption, $N^{+} \unlhd H^{+}$ and so for all $a^{+},b^{+}\in H^{+}$,
 \[(a^{+}\oplus b^{+})\oplus N^{+}=(a^{+}\oplus N^{+})\oplus b^{+}=a^{+}\oplus (N^{+} \oplus b^{+}).\]
This shows that
\[(a^{\epsilon}\otimes b^{\delta})\otimes N^{+}=(a^{\epsilon}\otimes N^{+})\otimes b^{\delta}=a^{\epsilon}\otimes (N^{+} \otimes b^{\delta}).\]
\item[(ii)] \textit{ $(\epsilon,\delta ) = (+,-) \ or \ (-,+)$}. By definition,
\begin{eqnarray*}
  (a^{+}\otimes b^{-})\otimes N^{+} &=& (a^{+} \oplus b^{+})^-\otimes N^{+}=((a^{+} \oplus b^{+})\oplus N^{+})^- \\
  (a^{+}\otimes N^{+})\otimes b^{-} &=& (a^{+}\oplus N^{+})\otimes b^{-}=((a^{+}\oplus N^{+})\oplus b^{+})^- \\
  a^{+}\otimes (N^{+}\otimes b^{-}) &=& a^{+}\otimes(N^{+}\oplus b^{+})^-=(a^{+}\oplus(N^{+}\oplus b^{+}))^- \\
  (a^{-}\otimes b^{+})\otimes N^{+} &=& (a^{+} \oplus b^{+})^-\otimes N^{+}=((a^{+} \oplus b^{+})\oplus N^{+})^- \\
  (a^{-}\otimes N^{+})\otimes b^{+} &=& (a^{+}\oplus N^{+})^-\otimes b^{+}=((a^{+}\oplus N^{+})\oplus b^{+})^- \\
  a^{-}\otimes (N^{+}\otimes b^{+}) &=& a^{-}\otimes(N^{+}\oplus b^{+})=(a^{+}\oplus(N^{+}\oplus b^{+}))^-
\end{eqnarray*}
By our assumption, $N^{+} \unlhd H^{+}$ and so for all $a^{+},b^{+}\in H^{+}$,
 \[(a^{+}\oplus b^{+})\oplus N^{+}=(a^{+}\oplus N^{+})\oplus b^{+}=a^{+}\oplus (N^{+} \oplus b^{+}).\]
Since $\phi$ is bijective,
\[((a^{+}\oplus b^{+})\oplus N^{+})^-=((a^{+}\oplus N^{+})\oplus b^{+})^-=(a^{+}\oplus (N^{+} \oplus b^{+}))^-.\]
This shows that
\[(a^{\epsilon}\otimes b^{\delta})\otimes N^{+}=(a^{\epsilon}\otimes N^{+})\otimes b^{\delta}=a^{\epsilon}\otimes (N^{+} \otimes b^{\delta}).\]
\end{enumerate}
This proves that  $N^{+}\unlhd G$, as desired.
\end{proof}

\section{Proof of the Main Result}

The aim of this section is to prove the main result of this paper. To do this, we assume that $M$ is a normal subgyrogroup of the gyrogroup $G$ introduced in Section 1.
By Theorem \ref{t2.1}, for each $a^{\epsilon}, b^{\delta} \in G=H^{+}\cup H^{-}$,
\begin{equation}\label{e1}
(a^{\epsilon}\otimes b^{\delta})\otimes M = (a^{\epsilon}\otimes M)\otimes b^{\delta} = a^{\epsilon}\otimes (M\otimes b^{\delta}).
\end{equation}
By Theorem \ref{t2.2}, one of the following conditions hold:

\begin{enumerate}
\item[(a)] $M \leq H^+$;

\item[(b)] there exists $N^+ \leq H^+$ and $L^- \subseteq H^-$ such that $M = N^+ \cup L^-$ and for all $x, y \in L^-$, $x \otimes y \in N^+$. Also, $N^+ \cap L^+ = \emptyset$ and $N^+ \cup L^+ \leq H^+$, where $L^+ = \varphi^{-1}(L^-)$;

\item[(c)] $M = N^+ \cup N^-$ such that $N^- = \varphi(N^+)$.
\end{enumerate}

Suppose that the condition (a) is satisfied. Then by considering $\delta = \varepsilon = +$ in Equation (\ref{e1}),
$M \unlhd H^+$. If condition (b) is satisfied, then $H^+ \cap M$ $=$ $(H^+ \cap N^+) \cup (H^+ \cap L^-)$ $=$ $H^+ \cap N^+$ $=$ $N^+$. Note that by our assumption, $M$ is normal in $G$ and by Corollary \ref{c2.3}(1), $H^+ \unlhd G$.
Hence, by \cite[Theorem 2.2]{3}, $N^+ \unlhd G$, which implies that $N^+ \unlhd H^+$. To complete this part, it is enough to prove that $N^+ \cup L^+ \unlhd H^+$, where $L^+ = \varphi^{-1}(L^-)$. By our assumption, $M=N^{+}\cup L^{-}\unlhd G$ and by Equation (1), $(a^{\epsilon}\otimes b^{\delta})\otimes (N^{+}\cup L^{-})=(a^{\epsilon}\otimes (N^{+}\cup L^{-}))\otimes b^{\delta}=a^{\epsilon}\otimes ((N^{+}\cup L^{-})\otimes b^{\delta})$, where  $a^{\epsilon},b^{\delta}$ are arbitrary elements of $G$. Therefore,

\begin{eqnarray}\label{e2}
\nonumber ((a^{\epsilon} \otimes b^{\delta})\otimes N^{+} )\cup ((a^{\epsilon}\otimes b^{\delta})\otimes L^{-}) &=&((a^{\epsilon}\otimes N^{+})\otimes b^{\delta} )\cup ((a^{\epsilon}\otimes L^{-})\otimes b^{\delta})\\
&=&(a^{\epsilon}\otimes(N^{+}\otimes b^{\delta})) \cup (a^{\epsilon}\otimes (L^{-}\otimes b^{\delta})).
\end{eqnarray}

We know that $N^+ \unlhd H^+$, and by Corollary \ref{c2.3}(2), $N^+ \unlhd G$. So for all  $a^{\epsilon},b^{\delta}\in G$,
\begin{equation}\label{e3}
  (a^{\epsilon}\otimes b^{\delta})\otimes N^{+}=(a^{\epsilon}\otimes N^{+})\otimes b^{\delta}=a^{\epsilon}\otimes (N^{+}\otimes b^{\delta})
\end{equation}
Since $N^+ \cap L^-=\emptyset$, by Equations (\ref{e2}) and (\ref{e3}),
\[(a^{\epsilon}\otimes b^{\delta})\otimes L^{-}=(a^{\epsilon}\otimes L^{-})\otimes b^{\delta}=a^{\epsilon}\otimes (L^{-}\otimes b^{\delta}).\]
In the previous equation, we set $\epsilon=\delta=+$, then
\[(a^{+}\otimes b^{+})\otimes L^{-}=(a^{+}\otimes L^{-})\otimes b^{+}=a^{+}\otimes (L^{-}\otimes b^{+}).\]
By definition, $((a^{+}\oplus b^{+})\oplus L^{+})^{-}= ((a^{+}\oplus L^{+})\oplus b^{+})^{-}
=(a^{+}\oplus( L^{+}\oplus b^{+}))^{-}.$ Since $\varphi$ is bijective, $(a^{+}\oplus b^{+})\oplus L^{+}= (a^{+}\oplus L^{+})\oplus b^{+} = a^{+}\oplus( L^{+}\oplus b^{+}).$
By the last equality and Equation (\ref{e3}), $((a^{+} \oplus b^{+})\oplus N)$ $\cup$ $((a^{+}\oplus b^{+})\oplus L^{+})$ $=$ $((a^{+}\oplus N)\oplus b^{+})$ $\cup$ $((a^{+}\oplus L^{+})\oplus b^{+})$ $=$ $(a^{+}\oplus(N\oplus b^{+}))$ $\cup$ $(a^{+}\oplus( L^{+}\oplus b^{+}))$ and hence $(a^{+}\oplus b^{+})\oplus (N^{+}\cup L^{+})=(a^{+}\oplus (N^{+}\cup L^{+}))\oplus b^{+}=a^{+}\oplus ((N^{+}\cup L^{+})\oplus b^{+}).$ Then $N^{+}\cup L^{+} \unlhd H^{+}$. This completes the proof of $(2)$. To prove $(3)$, we have to show that $N^{+}  \unlhd H^{+}$. By condition $(c)$, $H^{+}\cap M$ $=$ $(H^{+}\cap N^{+})$ $\cup$ $(H^{+}\cap N^{-})$ $=$ $H^{+}\cap N^{+}=N^{+}.$ By our assumption, $M\unlhd G$ and by Corollary \ref{c2.3},
$H^{+} \unlhd  G$. We now apply \cite[Theorem 2.2]{3} to deduce that $N^{+}\unlhd G$. Therefore,  $N^{+}\unlhd H^{+} $, as desired.

Conversely, we assume that $M$ satisfies one of the conditions $(1)$, $(2)$ or $(3)$ in Theorem \ref{t1.1}. It will be shown that $M \unlhd G$. To do this, the following three cases will be considered:

\begin{enumerate}
\item[(A)]  $M \unlhd H^{+}$.  By Corollary \ref{c2.3}(2), $M\unlhd G$.

\item[(B)]  \textit{There exists $N^+ \unlhd H^+$ and $L^- \subseteq H^-$ such that $M = N^+ \cup L^-$ and for all $x, y \in L^-$, $x \otimes y \in N^+$. Also, $N^+ \cap L^+ = \emptyset$ and $N^+ \cup L^+ \unlhd H^+$, where $L^+ = \varphi^{-1}(L^-)$}.
    Since $N^{+}, N^{+} \cup L^{+}\unlhd H^{+}$, by Corollary \ref{c2.3}(2), $N^{+}, N^{+} \cup L^{+}\unlhd G$. Then by Theorem \ref{t2.1}, for all $a^{\epsilon}, b^{\delta} \in G$,
\begin{equation}\label{e4}
(a^{\epsilon}\otimes b^{\delta})\otimes N^{+} = (a^{\epsilon}\otimes N^{+})\otimes b^{\delta} = a^{\epsilon}\otimes (N^{+}\otimes b^{\delta})
\end{equation}
and  $(a^{\epsilon}\otimes b^{\delta})\otimes (N^{+} \cup L^{+}) = (a^{\epsilon}\otimes (N^{+} \cup L^{+}))\otimes b^{\delta} = a^{\epsilon}\otimes ((N^{+} \cup L^{+})\otimes b^{\delta}).$ Hence, $((a^{\epsilon} \otimes b^{\delta})\otimes N^{+})$ $\cup$ $((a^{\epsilon}\otimes b^{\delta})\otimes L^{+})$ $=$ $((a^{\epsilon}\otimes N^{+})\otimes b^{\delta})$ $\cup$ $((a^{\epsilon}\otimes L^{+})\otimes b^{\delta})$ $=$ $(a^{\epsilon}\otimes(N^{+}\otimes b^{\delta}))$ $\cup$ $(a^{\epsilon}\otimes( L^{+}\otimes b^{\delta}))$. Since $N^+ \cap L^+ = \emptyset$, by Equation (\ref{e4}) and the last equality,
\begin{equation}\label{e5}
  (a^{\epsilon}\otimes b^{\delta})\otimes L^{+} = (a^{\epsilon}\otimes L^{+})\otimes b^{\delta} = a^{\epsilon}\otimes (L^{+}\otimes b^{\delta}).
\end{equation}
Since $\varphi$ is bijective,
\begin{equation}\label{e6}
  ((a^{\epsilon}\otimes b^{\delta})\otimes L^{+})^- = ((a^{\epsilon}\otimes L^{+})\otimes b^{\delta})^- = (a^{\epsilon}\otimes (L^{+}\otimes b^{\delta}))^-.
\end{equation}
Now, we claim that for all $a^{\epsilon},b^{\delta}\in G$, the following equality holds:
\begin{equation}\label{e7}
  (a^{\epsilon}\otimes b^{\delta})\otimes L^{-} = (a^{\epsilon}\otimes L^{-})\otimes b^{\delta} = a^{\epsilon}\otimes (L^{-}\otimes b^{\delta}).
\end{equation}
We consider two cases as follows:
\begin{enumerate}
\item[(B1)] $(\epsilon,\delta)=(+,+) \ or \ (-,-)$. By definition,
\begin{eqnarray*}
  (a^{+}\otimes b^{+})\otimes L^{-} &=& (a^{+} \oplus b^{+})\otimes L^{-}=((a^{+} \oplus b^{+})\oplus L^{+})^- \\
  (a^{+}\otimes L^{-})\otimes b^{+} &=& (a^{+}\oplus L^{+})^-\otimes b^{+}=((a^{+}\oplus L^{+})\oplus b^{+})^- \\
  a^{+}\otimes (L^{-}\otimes b^{+}) &=& a^{+}\otimes(L^{+}\oplus b^{+})^-=(a^{+}\oplus(L^{+}\oplus b^{+}))^- \\
  (a^{-}\otimes b^{-})\otimes L^{-} &=& (a^{+} \oplus b^{+})\otimes L^{-}=((a^{+} \oplus b^{+})\oplus L^{+})^- \\
  (a^{-}\otimes L^{-})\otimes b^{-} &=& (a^{+}\oplus L^{+})^-\otimes b^{+}=((a^{+}\oplus L^{+})\oplus b^{+})^- \\
  a^{-}\otimes (L^{-}\otimes b^{-}) &=& a^{-}\otimes(L^{+}\oplus b^{+})^-=(a^{+}\oplus(L^{+}\oplus b^{+}))^-.
\end{eqnarray*}
In this case, by Equation (\ref{e6}), our claim is true.
\item[(B2)] $(\epsilon,\delta) = (+,-) \ or \ (-,+)$. By definition,
\begin{eqnarray*}
  (a^{+}\otimes b^{-})\otimes L^{-} &=& (a^{+} \oplus b^{+})^-\otimes L^{-}=(a^{+} \oplus b^{+})\oplus L^{+} \\
  (a^{+}\otimes L^{-})\otimes b^{-} &=& (a^{+}\oplus L^{+})^-\otimes b^{-}=(a^{+}\oplus L^{+})\oplus b^{+} \\
  a^{+}\otimes (L^{-}\otimes b^{-}) &=& a^{+}\otimes(L^{+}\oplus b^{+})=a^{+}\oplus(L^{+}\oplus b^{+}) \\
  (a^{-}\otimes b^{+})\otimes L^{-} &=& (a^{+} \oplus b^{+})^-\otimes L^{-}=(a^{+} \oplus b^{+})\oplus L^{+} \\
  (a^{-}\otimes L^{-})\otimes b^{+} &=& (a^{+}\oplus L^{+})\otimes b^{+}=(a^{+}\oplus L^{+})\oplus b^{+} \\
  a^{-}\otimes (L^{-}\otimes b^{+}) &=& a^{-}\otimes(L^{+}\oplus b^{+})^-=a^{+}\oplus(L^{+}\oplus b^{+})
\end{eqnarray*}
Also, in this case, by Equation (\ref{e5}),  our claim is true.
\end{enumerate}
By Equations (\ref{e4}) and (\ref{e7}), we can see that $((a^{\epsilon} \otimes b^{\delta})\otimes N^{+})$ $\cup$ $((a^{\epsilon}\otimes b^{\delta})\otimes L^{-})$ $=$ $((a^{\epsilon}\otimes N^{+})\otimes b^{\delta})$ $\cup$ $((a^{\epsilon}\otimes L^{-})\otimes b^{\delta})$ $=$ $(a^{\epsilon}\otimes(N^{+}\otimes b^{\delta}))$ $\cup$ $(a^{\epsilon}\otimes( L^{-}\otimes b^{\delta}))$, which is equivalent to $((a^{\epsilon} \otimes b^{\delta})\otimes (N^{+} \cup  L^{-}))$ $=$ $((a^{\epsilon}\otimes (N^{+}\cup  L^{-} ))\otimes b^{\delta})$ $=$ $(a^{\epsilon}\otimes((N^{+} \cup  L^{-})\otimes b^{\delta}))$ or
$((a^{\epsilon} \otimes b^{\delta})\otimes M)$ $=$ $((a^{\epsilon}\otimes M)\otimes b^{\delta})$
$=$ $(a^{\epsilon}\otimes( M\otimes b^{\delta}))$.
This proves that $M \unlhd G$.


\item[(C)]  \textit{$M = N^+ \cup N^-$ such that $N^+ \unlhd H^+$ and $N^- = \varphi(N^+)$}. By Theorem \ref{t2.2}, $M \leq G$.
    Since $N^+ \unlhd H^+$, by Corollary \ref{c2.3}, $N^+ \unlhd G$. By Theorem \ref{t1.1}, for all $a^{\epsilon}, b^{\delta}\in G$,
\begin{equation}\label{e8}
  (a^{\epsilon}\otimes b^{\delta})\otimes N^+ = (a^{\epsilon}\otimes N^+)\otimes b^{\delta}=a^{\epsilon}\otimes (N^+\otimes b^{\delta}).
\end{equation}
 Since $\varphi$ is bijective,
  $((a^{\epsilon}\otimes b^{\delta})\otimes N^{+})^- = ((a^{\epsilon}\otimes N^{+})\otimes b^{\delta})^- = (a^{\epsilon}\otimes (N^{+}\otimes b^{\delta}))^-.$ Similar to the part B, we can show that $(a^{\epsilon}\otimes b^{\delta})\otimes N^{-}$ $=$ $(a^{\epsilon}\otimes N^{-})\otimes b^{\delta}$ $=$ $a^{\epsilon}\otimes (N^{-}\otimes b^{\delta})$. By the last equality and Equation (\ref{e8}), $(a^{\epsilon} \otimes b^{\delta})\otimes N^{+}$ $\cup$ $(a^{\epsilon}\otimes b^{\delta})\otimes N^{-}$ $=$ $(a^{\epsilon}\otimes N^{+})\otimes b^{\delta}$ $\cup$ $(a^{\epsilon}\otimes N^{-})\otimes b^{\delta}$ $=$ $a^{\epsilon}\otimes(N^{+}\otimes b^{\delta})$ $\cup$ $a^{\epsilon}\otimes( N^{-}\otimes b^{\delta})$ and so $(a^{\epsilon}\otimes b^{\delta})\otimes (N^{+} \cup N^{-})$ $=$ $(a^{\epsilon}\otimes (N^{+}\cup  N^{-}))\otimes b^{\delta}$ $=$ $a^{\epsilon}\otimes(( N^{+} \cup  N^{-})\otimes b^{\delta})$. We now apply our assumption to deduce that $(a^{\epsilon}\otimes b^{\delta})\otimes M$ $=$ $(a^{\epsilon}\otimes M)\otimes b^{\delta}$ $=$ $a^{\epsilon}\otimes(M\otimes b^{\delta})$, which means that  $M\unlhd G$.
\end{enumerate}

\section{Concluding Remarks}
In this paper, the normal subgyrogroups of a class of finite gyrogroups  is characterized. In this section, we check our main result by  a Gap code on three examples that is introduced in \cite{1}. Our Gap code is accessible from the authors upon request. Suppose that $K(1)$ is the gyrogroup such that its Cayley table is given in  Table \ref{T1}. We apply our method to construct the gyrogroups $K(2)$ from $K(1)$ and hence $|K(1)| = 8$  and $|K(2)| = 16$.

\begin{table}[H]
\centering
\caption{The Cayley Table of $K(1)$  such that $A=(4,5)(6,7)$.}\label{T1}
\vspace{-.3cm}
\begin{tabular}{c|cccccccc||c|cccccccc|}
$\oplus$ & \textbf 0 & \textbf 1 & \textbf 2 & \textbf 3 & \textbf 4 & \textbf 5 &\textbf  6 & \textbf 7 & $gyr_K$ &\textbf  0 & \textbf 1 &\textbf  2 &\textbf  3 &\textbf 4 &\textbf 5 &\textbf 6 &\textbf 7 \\\hline
\textbf 0 & 0 & 1 & 2 & 3 & 4 & 5 & 6 & 7 & \textbf 0 & $I$ & $I$ & $I$ & $I$ & $I$ & $I$ & $I$ & $I$\\
\textbf 1 & 1 & 0 & 3 & 2 & 5 & 4 & 7 & 6 &\textbf  1 & $I$ & $I$ & $I$ & $I$ & $I$ & $I$ & $I$ & $I$\\
\textbf 2&2&3&0&1&6&7&4&5 & \textbf 2& $I$ & $I$ & $I$ & $I$ & $A$ & $A$ & $A$ & $A$\\
\textbf 3&3&2&1&0&7&6&5&4 &\textbf  3& $I$ & $I$ & $I$ & $I$ & $A$ & $A$ & $A$ & $A$ \\
\textbf 4&4&5&6&7&0&1&2&3 & \textbf 4& $I$ & $I$ & $A$ & $A$ & $I$ & $I$ & $A$ & $A$\\
\textbf 5&5&4&7&6&1&0&3&2 &\textbf  5& $I$ & $I$ & $A$ & $A$ & $I$ & $I$ & $A$ & $A$\\
\textbf 6&6&7&4&5&3&2&1&0 & \textbf 6& $I$ & $I$ & $A$ & $A$ & $A$ & $A$ & $I$ & $I$\\
\textbf 7&7&6&5&4&2&3&0&1 &\textbf 7 & $I$ & $I$ & $A$ & $A$ & $A$ & $A$ & $I$ & $I$\\ \hline
\end{tabular}
\end{table}
\vspace{-.3cm}
\begin{table}[H]
\centering
\caption{The addition Table of $K(2)$.}
\vspace{-.3cm}
\begin{tabular}{ c|  c  c  c  c  c  c  c  c  c  c  c  c  c  c  c  c | }
            \hline
 $\oplus$ & $\textbf{0}$ & $\textbf{1}$ & $\textbf{2}$ & $\textbf{3}$ & $\textbf{4}$ & $\textbf{5}$ & $\textbf{6}$ & $\textbf{7}$ & $\textbf{8}$ & $\textbf{9}$ & $\textbf{10}$ & $\textbf{11}$ & $\textbf{12}$ & $\textbf{13}$ & $\textbf{14}$ & $\textbf{15}$ \\
\hline
        $\textbf{0}$ & $0$ & $1$ & $2$ & $3$ & $4$ & $5$ & $6$ & $7$ & $8$ & $9$ & $10$ & $11$ & $12$ & $13$ & $14$ & $15$ \\      
        $\textbf{1}$ &  $1$ & $0$ & $3$ & $2$ & $5$ & $4$ & $7$ & $6$ & $9$ & $8$ & $11$ & $10$ & $13$ & $12$ & $15$ & $14$ \\
         $\textbf{2}$ & $2$ & $3$ & $0$ & $1$ & $6$ & $7$ & $4$ & $5$ & $10$ & $11$ & $8$ & $9$ & $14$ & $15$ & $12$ & $13$   \\
         $\textbf{3}$ & $3$ & $2$ & $1$ & $0$ & $7$ & $6$ & $5$ & $4$ & $11$ & $10$ & $9$ & $8$ & $15$ & $14$ & $13$ & $12$  \\
         $\textbf{4}$ & $4$ & $5$ & $6$ & $7$ & $0$ & $1$ & $2$ & $3$ & $12$ & $13$ & $14$ & $15$ & $8$ & $9$ & $10$ & $11$  \\
         $\textbf{5}$ & $5$ & $4$ & $7$ & $6$ & $1$ & $0$&  $3$ & $2$ & $13$ & $12$ & $15$ & $14$ & $9$ & $8$ & $11$ & $10$  \\
         $\textbf{6}$ & $6$ & $7$ & $4$ & $5$ & $3$ & $2$ & $1$ & $0$ & $14$ & $15$ & $12$ & $13$ & $11$ & $10$ & $9$ & $8$  \\
         $\textbf{7}$ & $7$ & $6$ & $5$ & $4$ & $2$ & $3$ & $0$ & $1$ & $15$ & $14$ & $13$ & $12$ & $10$ & $11$ & $8$ & $9$ \\
         $\textbf{8}$ & $8$ & $9$ & $10$ & $11$ & $12$ & $13$ & $14$ & $15$ & $0$ & $1$ & $2$ & $3$ & $4$ & $5$ & $6$ & $7$  \\
         $\textbf{9}$ & $9$ & $8$ & $11$ & $10$ & $13$ & $12$ & $15$ & $14$ & $1$ & $0$ & $3$ & $2$ & $5$ & $4$ & $7$ & $6$  \\
         $\textbf{10}$ & $10$ & $11$ & $8$ & $9$ & $14$ & $15$ & $12$ & $13$ & $2$ & $3$ & $0$ & $1$ & $6$ & $7$ & $4$ & $5$ \\
         $\textbf{11}$ & $11$ & $10$ & $9$ & $8$ & $15$ & $14$ & $13$ & $12$ & $3$ & $2$ & $1$ & $0$ & $7$ & $6$ & $5$ & $4$  \\
         $\textbf{12}$ & $12$ & $13$& $14$ & $15$ & $8$ & $9$ & $10$ & $11$ & $4$ & $5$ & $6$ & $7$ & $0$ & $1$ & $2$ & $3$  \\
         $\textbf{13}$ & $13$ & $12$ & $15$ & $14$ & $9$ & $8$ & $11$ & $10$ & $5$ & $4$ & $7$ & $6$ & $1$ & $0$ & $3$ & $2$  \\
         $\textbf{14}$ & $14$ & $15$ & $12$ & $13$ & $11$ & $10$ & $9$ & $8$ & $6$ & $7$ & $4$ & $5$ & $3$ & $2$ & $1$ & $0$ \\
         $\textbf{15}$ & $15$ & $14$ & $13$ & $12$ & $10$ & $11$ & $8$ & $9$ & $7$ & $6$ & $5$ & $4$ & $2$ & $3$ & $0$ & $1$  \\
 \hline
 \end{tabular}
\end{table}
\begin{table}[htp]
\centering
\caption{The gyration table of $K(2)$ such that $A=(4,5)(6,7)(12,13)(14,15)$.}
\begin{tabular}{ c | c  c  c  c  c  c  c  c  c  c  c  c  c  c  c  c| }
 \hline
 $\textbf{gyr}$ & $\textbf{0}$ & $\textbf{1}$ & $\textbf{2}$ & $\textbf{3}$ & $\textbf{4}$ & $\textbf{5}$ & $\textbf{6}$ & $\textbf{7}$ & $\textbf{8}$ & $\textbf{9}$ & $\textbf{10}$ & $\textbf{11}$ & $\textbf{12}$ & $\textbf{13}$ & $\textbf{14}$ & $\textbf{15}$\\
 \hline
       $0$ & $I$ & $I$ & $I$ & $I$ & $I$ & $I$ & $I$ & $I$ & $I$ & $I$ & $I$ & $I$ & $I$ & $I$ & $I$ & $I$\\
       $\textbf{1}$ & $I$ & $I$ & $I$ & $I$ & $I$ & $I$ & $I$ & $I$ & $I$ & $I$ & $I$ & $I$ & $I$ & $I$ & $I$ & $I$\\
       $\textbf{2}$ & $I$ & $I$ & $I$ & $I$ & $A$ & $A$ & $A$ & $A$ & $I$ & $I$ & $I$ & $I$ & $A$ & $A$ & $A$ & $A$\\
       $\textbf{3}$ & $I$ & $I$ & $I$ & $I$ & $A$ & $A$ & $A$ & $A$ & $I$ & $I$ & $I$ & $I$ & $A$ & $A$ & $A$ & $A$\\
       $\textbf{4}$ & $I$ & $I$ & $A$ & $A$ & $I$ & $I$ & $A$ & $A$ & $I$ & $I$ & $A$ & $A$ & $I$ & $I$ & $A$ & $A$\\
       $\textbf{5}$  & $I$ & $I$ & $A$ & $A$ & $I$ & $I$ & $A$ & $A$ & $I$ & $I$ & $A$ & $A$ & $I$ & $I$ & $A$ & $A$\\
       $\textbf{6}$ & $I$ & $I$ & $A$ & $A$ & $A$ & $A$ & $I$ & $I$ & $I$ & $I$ & $A$ & $A$ & $A$ & $A$ & $I$ & $I$\\
       $\textbf{7}$ & $I$ & $I$ & $A$ & $A$ & $A$ & $A$ & $I$ & $I$ & $I$ & $I$ & $A$ & $A$ & $A$ & $A$ & $I$ & $I$\\
       $\textbf{8}$ & $I$ & $I$ & $I$ & $I$ & $I$ & $I$ & $I$ & $I$ & $I$ & $I$ & $I$ & $I$ & $I$ & $I$ & $I$ & $I$\\
       $\textbf{9}$ & $I$ & $I$ & $I$ & $I$ & $I$ & $I$ & $I$ & $I$ & $I$ & $I$ & $I$ & $I$ & $I$ & $I$ & $I$ & $I$\\
       $\textbf{10}$ & $I$ & $I$ & $I$ & $I$ & $A$ & $A$ & $A$ & $A$ & $I$ & $I$ & $I$ & $I$ & $A$ & $A$ & $A$ & $A$\\
       $\textbf{11}$ & $I$ & $I$ & $I$ & $I$ & $A$ & $A$ & $A$ & $A$ & $I$ & $I$ & $I$ & $I$ & $A$ & $A$ & $A$ & $A$\\
       $\textbf{12}$ & $I$ & $I$ & $A$ & $A$ & $I$ & $I$ & $A$ & $A$ & $I$ & $I$ & $A$ & $A$ & $I$ & $I$ & $A$ & $A$\\
       $\textbf{13}$ & $I$ & $I$ & $A$ & $A$ & $I$ & $I$ & $A$ & $A$ & $I$ & $I$ & $A$ & $A$ & $I$ & $I$ & $A$ & $A$\\
       $\textbf{14}$ & $I$ & $I$ & $A$ & $A$ & $A$ & $A$ & $I$ & $I$ & $I$ & $I$ & $A$ & $A$ & $A$ & $A$ & $I$ & $I$\\
       $\textbf{15}$ & $I$ & $I$ & $A$ & $A$ & $A$ & $A$ & $I$ & $I$ & $I$ & $I$ & $A$ & $A$ & $A$ & $A$ & $I$ & $I$\\
 \hline
\end{tabular}
\end{table}

\vskip 3mm

Our calculations show that all normal subgyrogroups of $K(1)$ and $K(2)$ are in Table \ref{T4}, respectively. In both cases, these normal subgyrogroups can be obtained from our man result. A nondegenerate gyrogroup is a gyrogroup that is not a group. In Table \ref{T4}, the subgyrogroups which is bolded are nondegenerate normal subgyrogroups.

\begin{table}[H]
\centering
\caption{The normal subgyrogroups of $K(1)$  and $K(2)$.}\label{T4}
\vspace{-.3cm}
\begin{tabular}{|c|c|}
\hline
$\textbf{Gyrogroup}$ & {\centering\textbf{The normal subgyrogroups}} \\\hline
\textbf{K(1)} & $\{0\}$, $\{0,1\}$, $\{0,1,2,3\}$, $\{0,1,4,5\}$, $\{0,1,6,7\}$,
              \textbf{\{0,1,2,3,4,5,6,7\}} \\ \hline
              & $\{0\}$, $\{0,1\}$, $\{0,1,2,3\}$, $\{0,1,4,5\}$, $\{0,1,6,7\}$,
              \textbf{\{0,1,2,3,4,5,6,7\}},\\
              & $\{0,8\}$, $\{0,9\}$, $\{0,1,8,9\}$, $\{0,1,10,11\}$, $\{0,1,12,13\}$, $\{ 0,1,14,15\}$, \\
\textbf{K(2)} &  $\{0,1,2,3,8,9,10,11\}$,\textbf{\{0,1,2,3,12,13,14,15\}},\\
              & $\{0,1,4,5,8,9,12,13 \}$, \textbf{\{0,1,4,5,10,11,14,15\}},\\
              & $\{ 0,1,6,7,8,9,14,15\}$, \textbf{\{0,1,6,7,10,11,12,13\}}, \textbf{\{0,1,2,3, \ldots, 13,14,15\}} \\ \hline
\end{tabular}
\end{table}

\vskip 3mm



\end{document}